\documentclass[authoryear]{paper}
\usepackage[T1]{fontenc}
\usepackage{geometry}
\usepackage{color}
\usepackage{array}
\usepackage{multirow}
\usepackage{amsmath}
\usepackage{graphicx}
\usepackage[authoryear]{natbib}

\makeatletter

\providecommand{\tabularnewline}{\\}

\@ifundefined{showcaptionsetup}{}{%
 \PassOptionsToPackage{caption=false}{subfig}}
\usepackage{subfig}
\makeatother

\begin{document}

\title{{\Large{}An optimisation approach for fuel treatment planning to
break the connectivity of high-risk regions}}

\author{\textcolor{black}{\large{}Ramya Rachmawati, Melih Ozlen, Karin J.
Reinke, John W. Hearne}}
\maketitle
\begin{abstract}
Uncontrolled wildfires can lead to loss of life and property and destruction
of natural resources. At the same time, fire plays a vital role in
restoring ecological balance in many ecosystems. Fuel management,
or treatment planning by way of planned burning, is an important tool
used in many countries where fire is a major ecosystem process. In
this paper, we propose an approach to reduce the spatial connectivity
of fuel hazards while still considering the ecological fire requirements
of the ecosystem. A mixed integer programming (MIP) model is formulated
in such a way that it breaks the connectivity of high-risk regions
as a means to reduce fuel hazards in the landscape. This multi-period
model tracks the age of each vegetation type and determines the optimal
time and locations to conduct fuel treatments. The minimum and maximum
Tolerable Fire Intervals (TFI), which define the ages at which certain
vegetation type can be treated for ecological reasons, are taken into
account by the model. Previous work has been limited to using single
vegetation types implemented within rectangular grids. In this paper,
we significantly extend previous work by modelling multiple vegetation
types implemented within a polygon-based network. Thereby a more realistic
representation of the landscape is achieved. An analysis of the proposed
approach was conducted for a fuel treatment area comprising 711 treatment
units in the Barwon-Otway district of Victoria, Australia. The solution
of the proposed model can be obtained for 20-year fuel treatment planning
within a reasonable computation time of eight hours.\end{abstract}

\section{Introduction}

Uncontrolled wildfires can result in the loss of life and economic
assets and the destruction of natural resources \citep{King2008421}.
Southern Australia, Mediterranean Europe and areas of the United States
are among the top regions in the world that are affected by frequent
wildfires \citep{Bradstock2012Wildfires}. Coupled with the proximity
of major cities to natural ecosystems prone to wildfire, the management
of fuel hazard becomes an important land management policy and planning
issue for the protection of human life and assets \citep{Collins201024}.
However, it is also acknowledged that fuel management for asset protection
cannot be done in isolation of the ecological requirements of the
ecosystem. Maintaining the ecological integrity of the landscape must
also be considered \citep{penman2011prescribed}.

Fuel management is a method to modify the structure and amount of
fuel. The methods include prescribed burning and mechanical clearing
\citep{King2008421,Loehle2004261}. Fuel management programs have
been extensively implemented in the USA \citep{Ager2010_FEM,Collins201024}
and Australia \citep{Boer2009132,McCaw2013217} in an effort to lessen
the risk posed by wildfire. The choice of fuel treatment location
plays a substantial role in conducting efficient fuel treatment scheduling
\citep{Collins201024}. Instead of randomly selecting the locations,
significantly better protection in a landscape could be provided by
a fuel treatment schedule that takes into account the relationships
between treatment units \citep{Schmidt20083170}. Research indicates
that it is important to choose where to conduct the fuel treatment
by considering spatial arrangement \citep{Rytwinski2010,Kim2009253,chung2015optimizing}.
The importance of landscape-level fuel treatment has been observed
in a number of studies. In wilderness regions in the United States,
a mosaic of varying fuel ages is formed as a result of free burning
fires. A particular arrangement of old and new treatment units has
been recognised to delay large wildfires in the following year \citep{Finney2007702}.
Research conducted in the Sierra Nevada forests of the United States
has shown that wildfire size can be modified by spatial fragmentation
of fuel \citep{vanwagtendonk1995large}. Prescribed burning has been
implemented in the eucalypt forests in south-western Australia over
the past 50 years. The connectivity of \textquoteleft old\textquoteright{}
untreated patches has been revealed to be the main aspect that contributes
to wildfire extent \citep{Boer2009132}.

Previous studies have proposed mathematically modelling fuel treatment
schedules methods for reducing fuel hazards. The studies had different
objective functions and took into account various considerations in
building up the models. \textcolor{black}{\citet{ferreira2011sto}}
propose a stochastic dynamic programming (SDP) approach to determine
the fuel treatment scheduling that produces the maximum expected discounted
net revenue while mitigating the risk of fire. The method was then
applied to a maritime pine forest in Leiria National Forest, Portugal.
They found that the approach was efficient and can efficiently help
integrating wildfire risk in stand management planning. \citet{garcia2011}
use the Hooke-Jeeves direct search method to determine the optimal
fuel treatment scheduling for reducing expected damage and increasing
the revenue to the same landscape, as that of \textcolor{black}{\citet{ferreira2011sto}}.
Their research shows that the fuel treatments improve productivity
as well as reduce the potential damage. \citet{Rachmawati2015} propose
a model that can lessen the risk of fire by reducing the total fuel
load but do not consider spatial properties or the spatial relationship
between the treatment units. \citet{wei2014schedule} propose a single-period
model to fragment high-risk patches by considering future fire spread
speeds and durations. \citet{Minas2014412} propose a model that breaks
the connectivity of high fuel units in the landscape to prevent the
fires spreading. The model proposed by \citet{Minas2014412} takes
into account vegetation dynamics in the landscape, but this is limited
to a simplistic grid representation of a single vegetation type per
treatment unit. In reality, a treatment unit may comprise a number
of patches with different vegetation type and age. In summary, most
of the models reviewed can be improved by taking into account multi-vegetation
types and ages within a treatment unit and using a polygon-based network
representation.

In this paper, we build upon previous work by incorporating multiple
vegetation types found in the landscape and within single treatment
units, and take into account the spatial connectivity or fragmentation
of \textquoteleft high-risk\textquoteright{} treatment units. We also
use a more realistic polygon-based network representation of the landscape
to better capture the spatial complexity of this problem rather than
a rectangular grid. Besides the negative impacts of wildfires, the
role of fire in ecology has been widely acknowledged. Fire is required
to maintain a healthy ecosystem and it also has a significant role
in habitat regeneration. Many vegetation species in fire-adapted ecosystems
need fire to reproduce. For instance, germination of seeds and successful
establishment of plants in the jarrah forests of Western Australia
is very rarely found without fire intervention \citep{burrows2003fire}.
More recently, \citet{Burrows20082394} argues that fuel management
is important to support biodiversity conservation as well as to reduce
the negative impact of wildfires. A recognition of vegetation dynamics
over time is crucial in the planning of fuel treatment \citep{Krivtsov20092915}.
In this proposed model, the ecological fire requirements of each vegetation
type can be described using the minimum and maximum Tolerable Fire
Intervals (TFI). The minimum TFI is the minimum time required between
two consecutive fire events at a location and is based on the time
to reach maturity of the sensitive species in the vegetation class.
The maximum TFI refers to the maximum time needed between two fire
events at a location that considers the fire interval required for
fire-adapted species rejuvenation \citep{Cheal2010}. In this paper,
we assume that treating of vegetation whose age is between these two
intervals will maintain species diversity and hence support the ecosystem\textquoteright s
health. Therefore, we select not to treat a treatment unit if the
age of vegetation growing in that location is under the minimum TFI.
In contrast, treatment units with vegetation over the maximum TFI
must be treated. In this paper, we assume that the high-risk threshold
age is between these two intervals. The objective of the model proposed
in this paper is to reduce the spatial connectivity of fuel hazards
while still considering the fire requirements of the ecosystem. The
question that then arises is when and where to conduct fuel treatment
to meet this objective, that can be solved for spatially complex landscapes
with long planning horizons? 

A Mixed Integer Programming (MIP) model is proposed for multi-period
fuel treatment scheduling. The model tracks the vegetation age in
each treatment unit yearly for both treated and untreated areas. The
model is then applied to a real landscape in southern Australia that
comprises different shapes and sizes of treatment units.

\section{Problem formulation}

In this section, we explain the terms \textquoteleft treatment unit\textquoteright{}
and \textquoteleft patch\textquoteright{} that we use to formulate
the problem. The candidate locations for fuel treatment are represented
by treatment units. A treatment unit comprises multiple patches. Each
vegetation type growing in a treatment unit is represented by a patch
and within each patch all the vegetation is of the same age. The data
in each patch includes area, vegetation type and age. Patches within
a single treatment unit may have different vegetation type and age,
defining a \textquoteleft multi-vegetation treatment unit\textquoteright . 

Each vegetation type has a \textquoteleft high risk\textquoteright{}
age threshold. For example, grass and bush are considered to be high
risk when they reach four and seven years old, respectively. Since
we know the vegetation type and age in each patch, we then know whether
a patch is a high-risk patch or not at any given time. In order to
disconnect the high-risk treatment units in a landscape, we need a
method to determine whether a treatment unit is a high-risk treatment
unit or not. In this paper, we assume that if ignitions occur, the
fires will likely spread through connected high-risk treatment units.
From this, we believe that if we can disconnect high-risk treatments
unit as much as possible, the possibility of catastrophic fires can
be reduced. 

Each treatment unit selected for fuel treatment should not violate
the ecological requirements. Each vegetation type has its specific
minimum and maximum TFI. We assume that a healthier ecosystem can
be maintained when the fuel treatment is conducted when the vegetation
age is between the minimum and the maximum TFI.

\section{Model formulation{\normalsize{}\label{sec:Model-formulation}}}

The model is formulated to determine when and where to conduct the
fuel treatment each year to break the connectivity of high-risk treatment
units and to meet the ecological requirements. We consider a landscape
divided into treatment units where each treatment unit might consist
of multiple patches. The following mixed integer programming model
is formulated. 

\textcolor{black}{\smallskip{}
}

\textcolor{black}{Sets:}

\textcolor{black}{\smallskip{}
}

\textcolor{black}{$C$ is the set of all treatment units in the landscape}

\textcolor{black}{$\Psi\subset C$ is the set of treatment units where
fuel treatment is not permitted}

\textcolor{black}{$\varLambda\subset C$ is the set of treatment units
where fuel treatment is permitted (where $\varLambda=C-\Psi$)}

\textcolor{black}{$P_{i}$ is the set of patches in treatment unit
$i$}

\textcolor{black}{$\varPhi_{i}$ is the set of treatment units connected
to treatment unit $i$}

\textcolor{black}{$T$ is the planning horizon}

\textcolor{black}{\smallskip{}
}

\textcolor{black}{Indices:}

\textcolor{black}{\smallskip{}
}

\textcolor{black}{$p$ = patch}

\textcolor{black}{$i$ = treatment unit}

\textcolor{black}{$t$ = period, $t$ = 0, 1, 2, \ldots}\textcolor{black}{\emph{T}}

\textcolor{black}{\smallskip{}
}

\textcolor{black}{Parameters:}

\textcolor{black}{\smallskip{}
}

\textcolor{black}{$w_{i,j}$ = relative importance (weight) of connectivity
of treatment units }\textcolor{black}{\emph{i }}\textcolor{black}{and
}\textcolor{black}{\emph{j}}

\textcolor{black}{$a_{p}$ = initial vegetation age in patch }\textcolor{black}{\emph{p}}

$Area_{p}$= area of patch {\textcolor{black}{\emph{p}}}

$\rho$ = treatment level (in percentage), i.e. the maximum proportion
of the total area 

\hspace{1cm}that fuel treatment is permitted in a landscape selected
for treatment

\textcolor{black}{$R$ = the total area of }{\textcolor{black}{treatment
units in the landscape where fuel treatment is permitted}}

$c_{i}$= area of treatment unit $i$

\textcolor{black}{$d_{p}$ = high-risk age threshold for patch }\textcolor{black}{\emph{p,}}\textcolor{black}{{}
based upon the vegetation type growing }

\hspace{1cm}\textcolor{black}{in that patch}

\textcolor{black}{$maxTFI_{p}$= maximum tolerable fire interval (TFI)
of vegetation type growing in patch }{\textcolor{black}{$p$}}

\textcolor{black}{$minTFI_{p}$ = minimum TFI of vegetation type growing
in patch }{\textcolor{black}{$p$}}

$H$ = the threshold for the area proportion of the high-risk patches
in a treatment unit 

\hspace{1cm}to be a high-risk treatment units

\textcolor{black}{\smallskip{}
}

\textcolor{black}{Decision variables:}

\textcolor{black}{$A_{p,t}$ = vegetation age in patch }\textcolor{black}{\emph{p}}\textcolor{black}{{}
at time }\textcolor{black}{\emph{t}}\textcolor{black}{{} }

$x_{i,t}=\begin{cases}
1 & \mbox{if treatment unit \emph{i} is treated in time period \emph{t}}\\
0 & \mbox{otherwise}
\end{cases}$

\textcolor{black}{\smallskip{}
}

$Riskpatch_{p,t}=\begin{cases}
1 & \mbox{if patch \emph{p} is classified as high-risk patch in time period \emph{t}}\\
0 & \mbox{otherwise}
\end{cases}$

\textcolor{black}{\smallskip{}
}

$Risk_{i,t}=\begin{cases}
1 & \mbox{if treatment unit \emph{i} is classified as high-risk treatment unit in time period \emph{t}}\\
0 & \mbox{otherwise}
\end{cases}$

\vspace{0.5cm}

$RiskConn_{i,j,t}=\begin{cases}
1 & \mbox{if connected treatment units \emph{i} and \emph{j} are both high-risk treatment units in time period \emph{t}}\\
0 & \mbox{otherwise}
\end{cases}$

\textcolor{black}{\smallskip{}
}

$Old_{p,t}=\begin{cases}
1 & \textrm{if patch \emph{p }is classified as\text{ \textquoteleft}old\text{\textquoteright\ }(over-the-maximum-TFI)}\\
 & \textrm{patch in time period \emph{t}}\\
0 & \textrm{otherwise}
\end{cases}$

\textcolor{black}{\smallskip{}
}

$Young_{p,t}=\begin{cases}
1 & \textrm{if patch \emph{p }is classified as\text{ \textquoteleft}young\text{\textquoteright\ }(under-the-minimum-TFI)}\\
 & \textrm{patch in time period \emph{t}}\\
0 & \textrm{otherwise}
\end{cases}$

\textcolor{black}{\smallskip{}
}

\textcolor{black}{\smallskip{}
}

\textcolor{black}{Minimise the weighted connectivity of high-risk
treatment units}

\textcolor{black}{
\begin{equation}
z=\underset{t=1}{\overset{T}{\sum}}\underset{i\in C}{\sum}\underset{j\in\varPhi_{i},i<j}{\sum}w_{i,j}RiskConn_{i,j,t}\label{eq:objective function 1}
\end{equation}
}

\textcolor{black}{subject to}

\textcolor{black}{
\begin{equation}
\underset{i}{\sum}c_{i}x_{i,t}\leq\rho R,\;t=1\ldots T,\forall i\in\varLambda\label{eq:budget or level treatment}
\end{equation}
\begin{equation}
A_{p,0}=a_{p},\;\forall p\label{eq:initial age}
\end{equation}
}

\textcolor{black}{
\begin{equation}
A_{p,t}=A_{p,t-1}+1,\;\forall p\in P_{i},t=1\ldots T,\forall i\in\varPsi\label{eq:veg age next period will be added by one for trmnt unit where fuel tr is not permitted}
\end{equation}
}

\textcolor{black}{
\begin{equation}
A_{p,t}\geq A_{p,t-1}+1-M_{1}x_{i,t},\;\forall p\in P_{i},t=1\ldots T,\forall i\in\varLambda\label{eq:veg age next period will be added by one if xit 0}
\end{equation}
}

\textcolor{black}{
\begin{equation}
A_{p,t}\leq M_{2}(1-x_{i,t}),\;\forall p\in P_{i},t=1\ldots T,\forall i\in\varLambda\label{eq:veg age is reset to zero if xit 1 (burned)}
\end{equation}
}

\textcolor{black}{
\begin{equation}
A_{p,t}\leq A_{p,t-1}+1,\;\forall p\in P_{i},t=1\ldots T,\forall i\in\varLambda\label{eq:force_veg age will be added by one if not burned}
\end{equation}
}

\textcolor{black}{
\begin{equation}
A_{p,t}-d_{p}\leq M_{3}Riskpatch_{p,t}-1,\;\forall p\in P_{i},t=1\ldots T,\forall i\in C\label{eq:O_pt=00003D1 if high-risk}
\end{equation}
}

\textcolor{black}{
\begin{equation}
\underset{p\in P_{i}}{\overset{}{\sum}}Area_{p}Riskpatch_{p,t}-H\underset{p\in P_{i}}{\overset{}{\sum}}Area_{p}\leq M_{4}Risk_{i,t},\;t=1\ldots T,\forall p\in P_{i},\forall i\in C\label{eq:a standard for a treatment unit to be danger tr unit, H can be 50 percent}
\end{equation}
}

\textcolor{black}{
\begin{equation}
Risk_{i,t}+Risk_{j,t}-RiskConn_{i,j,t}\leq1,\;t=1\ldots T,\forall j\in\varPhi_{i},i<j,\forall i\in C\label{eq:connectivity of treatment unit}
\end{equation}
}

\textcolor{black}{
\begin{equation}
A_{p,t}-maxTFI_{p}\leq M_{5}Old_{p,t}-1,\;\forall p\in P_{i},t=0\ldots T-1,\forall i\in\varLambda\label{eq: tr unit i is classified 'Old'_TFI constraint}
\end{equation}
}

\textcolor{black}{
\begin{equation}
A_{p,t}\geq maxTFI_{p}Old_{p,t},\;\forall p\in P_{i},t=0\ldots T-1,\forall i\in\varLambda\label{eq: tr unit i is classified 'Old'_TFI constraint-1}
\end{equation}
}

\textcolor{black}{
\begin{equation}
A_{p,t}+M_{6}Young_{p,t}\geq minTFI_{p},\;\forall p\in P_{i},t=0\ldots T-1,\forall i\in\varLambda\label{eq: tr unit i is classified 'Old'_TFI constraint-1-1}
\end{equation}
}

\textcolor{black}{
\begin{equation}
A_{p,t}-M_{7}(1-Young_{p,t})\leq minTFI_{p}-1,\;\forall p\in P_{i},t=0\ldots T-1,\forall i\in\varLambda\label{eq: tr unit i is classified 'Old'_TFI constraint-1-1-1}
\end{equation}
}

\textcolor{black}{
\begin{equation}
Young{}_{p,t-1}\leq1-x{}_{i,t},\;t=1\ldots T,\forall i\in\varLambda\label{eq:if young 1 then x 0}
\end{equation}
}

\textcolor{black}{
\begin{equation}
\underset{p\in P_{i}}{\sum}Old{}_{p,t-1}-\mid V_{i}\mid\underset{p\in P_{i}}{\sum}Young{}_{p,t-1}\leq\mid V_{i}\mid x_{i,t},\;t=1\ldots T,\forall i\in\varLambda\label{eq:do not burn if there are any patches are young}
\end{equation}
}

\textcolor{black}{
\begin{equation}
x_{i,t},\;Riskpatch_{p,t},\;Risk_{i,t},\;RiskConn_{i,j,t},\;Young_{p,t},\;Old_{p,t}\in\{0,1\}\label{eq:binary vars}
\end{equation}
}

\textcolor{black}{The objective function \eqref{eq:objective function 1}
minimises the weighted connectivity of high-risk treatment units in
a landscape throughout a planning horizon. }

\textcolor{black}{Constraint \eqref{eq:budget or level treatment}}{
specifies that the total area selected for fuel treatment annually
is not more than the area allotted (target) each year for fuel treatment
(in hectares).}

\textcolor{black}{Constraint \eqref{eq:initial age} sets the initial
vegetation age in a patch. Constraint \eqref{eq:veg age next period will be added by one for trmnt unit where fuel tr is not permitted}
to \eqref{eq:veg age is reset to zero if xit 1 (burned)} track the
vegetation age of each patch. Constraint \eqref{eq:veg age next period will be added by one for trmnt unit where fuel tr is not permitted}
relates to the set of treatment units where fuel treatment is not
permitted. Constraint \eqref{eq:veg age next period will be added by one if xit 0}
and \eqref{eq:veg age is reset to zero if xit 1 (burned)} indicate
that when $x_{i,t}=0$, the vegetation in that area will continue
growing until the following period, and the age will be incremented
by one. Whereas if $x_{i,t}=1$, the vegetation age will reset to
zero. Constraint \eqref{eq:force_veg age will be added by one if not burned}
increments vegetation age by exactly one year if the treatment unit
is not treated.}

Constraint \eqref{eq:O_pt=00003D1 if high-risk} uses binary variable
\textcolor{black}{$Riskpatch_{p,t}$ to classify a patch to be a high-risk
patch if the vegetation age in that patch reaches or exceeds a threshold
value}, thus each patch has its own age threshold. Then, within a
single treatment unit, we can compare the area of over-the-threshold
patch. Here, we define a treatment unit is a high-risk treatment unit
if the proportion of the over the threshold area is greater than a
certain proportion of the total treatable area of the treatment unit.
Constraint \eqref{eq:a standard for a treatment unit to be danger tr unit, H can be 50 percent}
represents this requirement. In constraint \eqref{eq:connectivity of treatment unit},
$RiskConn_{i,j,t}$ takes the value one if connected treatment units
\emph{i} and \emph{j} are both classified as high-risk treatment units
in time period \emph{t} .

\textcolor{black}{Constraints \eqref{eq: tr unit i is classified 'Old'_TFI constraint}
to \eqref{eq: tr unit i is classified 'Old'_TFI constraint-1-1-1}
classify a patch to be an }{\textquoteleft old\textquoteright{}}\textcolor{black}{{}
or a }{\textquoteleft young\textquoteright{}}\textcolor{black}{{}
patch based on TFI values. Constraint \eqref{eq:if young 1 then x 0}
ensures that the treatment units containing young patches cannot be
treated. Constraint \eqref{eq:do not burn if there are any patches are young}
states that if there is at least one patch within a treatment unit
is }{\textquoteleft old\textquoteright{}}\textcolor{black}{{}
and no young patch, then the treatment unit must be treated. Here,
}{$|V_{i}|$}\textcolor{black}{{} represents
the number of patches in treatment unit $i$. This constraint avoids
a deadlock that may occur when a treatment unit consists of a young
and an old patch at the same time. In this study, we break the deadlock
in favour of young patch.}

\textcolor{black}{Constraints \eqref{eq:binary vars} ensures that
the decision variables take binary values.}

\subsection{Model improvements}

The solution time can be improved by reducing the number of variables.
As discussed earlier, the initial age of each vegetation type in each
treatment unit is given. We also assume that the age of vegetation
type growing in the treatment units where fuel treatment is not permitted
should always be incremented by one. For this reason, we no longer
need constraint \eqref{eq:veg age next period will be added by one for trmnt unit where fuel tr is not permitted}
to track the vegetation in the area. The time for the vegetation type
to reach the high-risk age threshold can be determined. \textcolor{black}{And
because we assume that we cannot treat the treatment units, once the
vegetation type hits the threshold it will remain high risk. Therefore,
within a planning horizon we can determine whether a treatment unit
is high risk or not. }

Decision variables \textcolor{black}{$A_{p,t}$} and \textcolor{black}{$Riskpatch_{p,t}$}
for the treatment units where fuel treatment is not permitted can
be omitted, and regarded as parameters instead. This results in a
faster solution time. 

We can rewrite our model as follows. Constraint \eqref{eq:veg age next period will be added by one for trmnt unit where fuel tr is not permitted}
is excluded, because at any given time the age of vegetation growing
in the treatment units where fuel treatment is not permitted is known.
Constraints \eqref{eq:O_pt=00003D1 if high-risk} and \eqref{eq:a standard for a treatment unit to be danger tr unit, H can be 50 percent}
are only defined for treatable treatment units. All other constraints
remain the same. However, we introduce these two constraints to the
model for the treatment units where fuel treatments are not permitted:

\bigskip{}

\textcolor{black}{
\begin{equation}
Risk_{i,t}=0,\:\forall t\:\textrm{when}\:\theta\leq0,\forall i\in\Psi\label{eq:no high risk}
\end{equation}
}

\textcolor{black}{
\begin{equation}
Risk_{i,t}=1,\:\forall t\:\textrm{when}\:\theta>0,\forall i\in\Psi\label{eq:high risk}
\end{equation}
}

where $\theta=\underset{p\in P_{i}}{\overset{}{\sum}}Area_{p}Riskpatch_{p,t}-H\underset{p\in P_{i}}{\overset{}{\sum}}Area_{p}$

\bigskip{}

In constraint \eqref{eq:no high risk}, value 0 is assigned to \textcolor{black}{$Risk_{i,t}$}
if less than a certain proportion of the total treatable area of the
treatment unit is high risk at time \emph{t}. And in constraint \eqref{eq:high risk}
value 1 is assigned to \textcolor{black}{$Risk_{i,t}$} if more than
a certain proportion of the total treatable area of the treatment
unit is high risk at time \emph{t}.

\section{\textcolor{black}{Implementation of the new approach}}

Initially, it may not be possible to treat all treatment units according
to the maximum TFI value because of the annual limit, $\rho$. This
maximum TFI requirement may lead to the infeasibility of the initial
problem. In order to bring the system under control and to avoid the
initial infeasibility, we propose a preliminary stage, namely Phase
1. From the initial data, we can identify treatment units containing
an old patch or would potentially be containing an old patch in the
following year and have no young patches. We are trying to eliminate
the treatment units containing old patches to ensure feasibility.
In this phase, we exclude the TFI constraints, which are constraints
\eqref{eq: tr unit i is classified 'Old'_TFI constraint} to \eqref{eq:do not burn if there are any patches are young}.
We run the model without enforcing the constraint ensuring treatment
of old patches for some years, and modify the objective function as
follows:

\textcolor{black}{maximise}

\textcolor{black}{
\begin{equation}
z=\underset{t=1}{\overset{N}{\sum}}\underset{i\in\varTheta}{\sum}c_{i}x_{i,t}-\underset{t=1}{\overset{N}{\sum}}\underset{i\in\varTheta}{\sum}\underset{j\in\varPhi_{i},i<j}{\sum}\varepsilon_{i}RiskConn_{i,j,t}\label{eq:phase1-obj-1}
\end{equation}
}

\textcolor{black}{\smallskip{}
}

where $\varTheta$ is the set of treatment units that contains an
old patch or potentially contains an old patch in the following year
and no young patch.{\textcolor{black}{{} $\varepsilon_{i}$
is a relatively small number ($\varepsilon_{i}\ll c_{i}$) representing
the weight of connectivity of treatment unit }}$i$. \emph{N} is the
planning horizon.{\textcolor{black}{{} }}

The objective is to maximise the area treated and to minimise the
weighted connectivity of the treatment units in a landscape for a
number of years ahead. The planning horizon (\emph{N}) increased incrementally
until the initial problem is feasible.

For the landscape that comprises mostly old treatment units, the solution
from this phase becomes the input for Phase 2. In Phase 2, the model
presented in Section \ref{sec:Model-formulation} is run.

\section{Model demonstration}

\begin{figure}
\protect\caption{A landscape for the model demonstration (29 treatment units)}

\subfloat[Map of the landscape\label{fig:Map-of-the29}]{\medskip{}

\protect\centering{}\protect\includegraphics[scale=0.3]{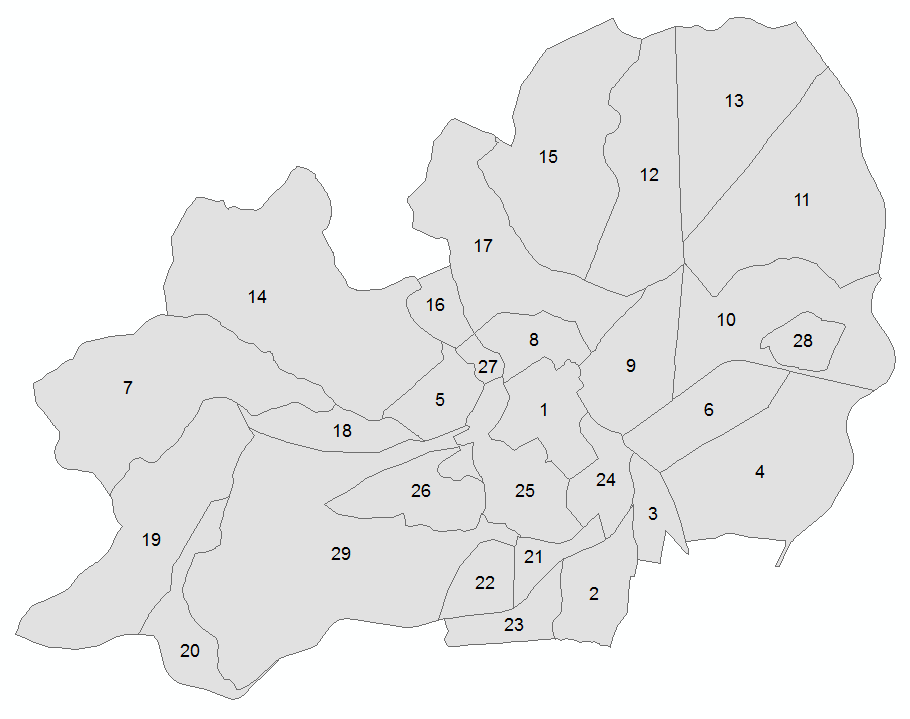}\protect}\subfloat[The neighbourhood graph of the landscape\label{fig:The-neighbourhood-graph-29}]{\medskip{}

\protect\centering{}\protect\includegraphics[scale=0.3]{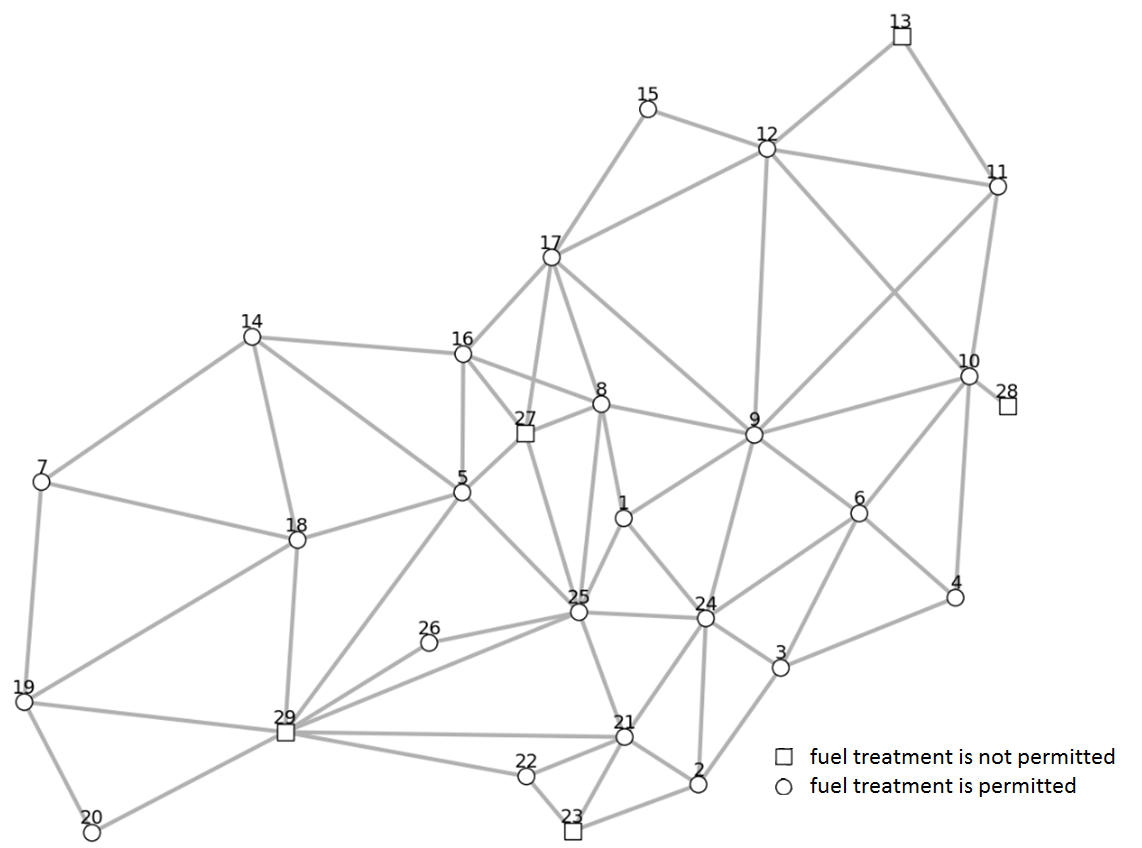}\protect}
\end{figure}

For the model demonstration, consider a test landscape comprising
29 treatment units that are a subset of the case study in the Barwon-Otway
district of Victoria, Australia. Figure \ref{fig:Map-of-the29} represents
the map of the landscape and Figure \ref{fig:The-neighbourhood-graph-29}
illustrates the graph representing the neighbourhood of each treatment
unit. We assume that two treatment units are neighbouring if they
have common boundaries. Table \ref{tab:Ecological-Vegetation-Class-29}
represents data for each Ecological Vegetation Class (EVC) and the
associated threshold age, the minimum and the maximum TFI for this
test landscape. The data regarding the area of the treatment units,
vegetation type (EVC) and age can be seen in Table \ref{tab:29-treatment-units}. 

\begin{table}
\protect\caption{Ecological Vegetation Class (EVC) and the associated threshold age,
the minimum and the maximum TFI for the test landscape\label{tab:Ecological-Vegetation-Class-29}}
\medskip{}

\centering{}%
\begin{tabular}{|c|c|c|c|}
\hline 
EVC code & min TFI (year) & max TFI (year) & threshold (year)\tabularnewline
\hline 
1 & 3 & 10 & 5\tabularnewline
\hline 
3 & 4 & 15 & 7\tabularnewline
\hline 
6 & 7 & 20 & 10\tabularnewline
\hline 
\end{tabular}
\end{table}

\begin{table}
\protect\caption{29 treatment units data containing vegetation type, extent and age\label{tab:29-treatment-units}}

\centering{}%
\begin{tabular}{|>{\centering}p{2cm}|>{\centering}p{0.7cm}|>{\centering}p{0.7cm}|>{\centering}p{1cm}|c|>{\centering}p{0.7cm}|>{\centering}p{0.7cm}|>{\centering}p{0.7cm}|>{\centering}p{0.7cm}|c|>{\centering}p{0.7cm}|>{\centering}p{0.7cm}|>{\centering}p{0.7cm}|>{\centering}p{0.7cm}|}
\cline{1-4} \cline{6-9} \cline{11-14} 
{\small{}Treatment unit ID} & {\small{}EVC code} & {\small{}area (ha)} & {\small{}age (years)} &  & {\small{}7} & {\small{}6} & {\small{}34} & {\small{}11} &  & {\small{}18} & {\small{}6} & {\small{}10} & {\small{}9}\tabularnewline
\cline{1-4} \cline{6-9} \cline{11-14} 
{\small{}1} & {\small{}1} & {\small{}10} & {\small{}6} &  & {\small{}8} & {\small{}1} & {\small{}19} & {\small{}5} &  & {\small{}19} & {\small{}1} & {\small{}50} & {\small{}5}\tabularnewline
\cline{1-4} \cline{6-9} \cline{11-14} 
{\small{}1} & {\small{}3} & {\small{}8} & {\small{}7} &  & {\small{}8} & {\small{}3} & {\small{}12} & {\small{}7} &  & {\small{}19} & {\small{}3} & {\small{}37} & {\small{}5}\tabularnewline
\cline{1-4} \cline{6-9} \cline{11-14} 
{\small{}1} & {\small{}6} & {\small{}14} & {\small{}11} &  & {\small{}9} & {\small{}1} & {\small{}46} & {\small{}4} &  & {\small{}20} & {\small{}1} & {\small{}10} & {\small{}1}\tabularnewline
\cline{1-4} \cline{6-9} \cline{11-14} 
{\small{}2} & {\small{}1} & {\small{}10} & {\small{}5} &  & {\small{}10} & {\small{}1} & {\small{}78} & {\small{}6} &  & {\small{}20} & {\small{}3} & {\small{}6} & {\small{}2}\tabularnewline
\cline{1-4} \cline{6-9} \cline{11-14} 
{\small{}2} & {\small{}3} & {\small{}21} & {\small{}8} &  & {\small{}11} & {\small{}1} & {\small{}30} & {\small{}4} &  & {\small{}20} & {\small{}6} & {\small{}14} & {\small{}10}\tabularnewline
\cline{1-4} \cline{6-9} \cline{11-14} 
{\small{}3} & {\small{}1} & {\small{}4} & {\small{}1} &  & {\small{}11} & {\small{}3} & {\small{}50} & {\small{}8} &  & {\small{}21} & {\small{}1} & {\small{}5} & {\small{}1}\tabularnewline
\cline{1-4} \cline{6-9} \cline{11-14} 
{\small{}3} & {\small{}3} & {\small{}5} & {\small{}1} &  & {\small{}11} & {\small{}6} & {\small{}30} & {\small{}12} &  & {\small{}21} & {\small{}3} & {\small{}8} & {\small{}1}\tabularnewline
\cline{1-4} \cline{6-9} \cline{11-14} 
{\small{}3} & {\small{}6} & {\small{}7} & {\small{}1} &  & {\small{}12} & {\small{}1} & {\small{}40} & {\small{}5} &  & {\small{}22} & {\small{}3} & {\small{}19} & {\small{}7}\tabularnewline
\cline{1-4} \cline{6-9} \cline{11-14} 
{\small{}4} & {\small{}1} & {\small{}40} & {\small{}5} &  & {\small{}12} & {\small{}3} & {\small{}34} & {\small{}7} &  & {\small{}23} & {\small{}6} & {\small{}20} & {\small{}11}\tabularnewline
\cline{1-4} \cline{6-9} \cline{11-14} 
{\small{}4} & {\small{}3} & {\small{}30} & {\small{}6} &  & {\small{}13} & {\small{}6} & {\small{}84} & {\small{}11} &  & {\small{}24} & {\small{}6} & {\small{}22} & {\small{}10}\tabularnewline
\cline{1-4} \cline{6-9} \cline{11-14} 
{\small{}4} & {\small{}6} & {\small{}24} & {\small{}10} &  & {\small{}14} & {\small{}3} & {\small{}80} & {\small{}7} &  & {\small{}25} & {\small{}1} & {\small{}42} & {\small{}1}\tabularnewline
\cline{1-4} \cline{6-9} \cline{11-14} 
{\small{}5} & {\small{}1} & {\small{}8} & {\small{}1} &  & {\small{}14} & {\small{}6} & {\small{}76} & {\small{}11} &  & {\small{}26} & {\small{}3} & {\small{}33} & {\small{}7}\tabularnewline
\cline{1-4} \cline{6-9} \cline{11-14} 
{\small{}5} & {\small{}3} & {\small{}10} & {\small{}1} &  & {\small{}15} & {\small{}6} & {\small{}103} & {\small{}12} &  & {\small{}27} & {\small{}3} & {\small{}6} & {\small{}6}\tabularnewline
\cline{1-4} \cline{6-9} \cline{11-14} 
{\small{}5} & {\small{}6} & {\small{}4} & {\small{}1} &  & {\small{}16} & {\small{}3} & {\small{}14} & {\small{}5} &  & {\small{}28} & {\small{}1} & {\small{}14} & {\small{}5}\tabularnewline
\cline{1-4} \cline{6-9} \cline{11-14} 
{\small{}6} & {\small{}1} & {\small{}18} & {\small{}1} &  & {\small{}17} & {\small{}1} & {\small{}50} & {\small{}5} &  & {\small{}29} & {\small{}1} & {\small{}100} & {\small{}5}\tabularnewline
\cline{1-4} \cline{6-9} \cline{11-14} 
{\small{}6} & {\small{}3} & {\small{}20} & {\small{}1} &  & {\small{}17} & {\small{}3} & {\small{}32} & {\small{}6} &  & {\small{}29} & {\small{}3} & {\small{}50} & {\small{}6}\tabularnewline
\cline{1-4} \cline{6-9} \cline{11-14} 
{\small{}7} & {\small{}3} & {\small{}80} & {\small{}8} &  & {\small{}18} & {\small{}3} & {\small{}14} & {\small{}5} &  & {\small{}29} & {\small{}6} & {\small{}41} & {\small{}9}\tabularnewline
\cline{1-4} \cline{6-9} \cline{11-14} 
\end{tabular}
\end{table}

\begin{figure}
\noindent \begin{centering}
\protect\caption{A network represents the fuel treatment schedule for the test landscape\label{fig:Fuel-treatment-schedule-29}}

\par\end{centering}

\medskip{}

\centering{}\includegraphics[scale=0.22]{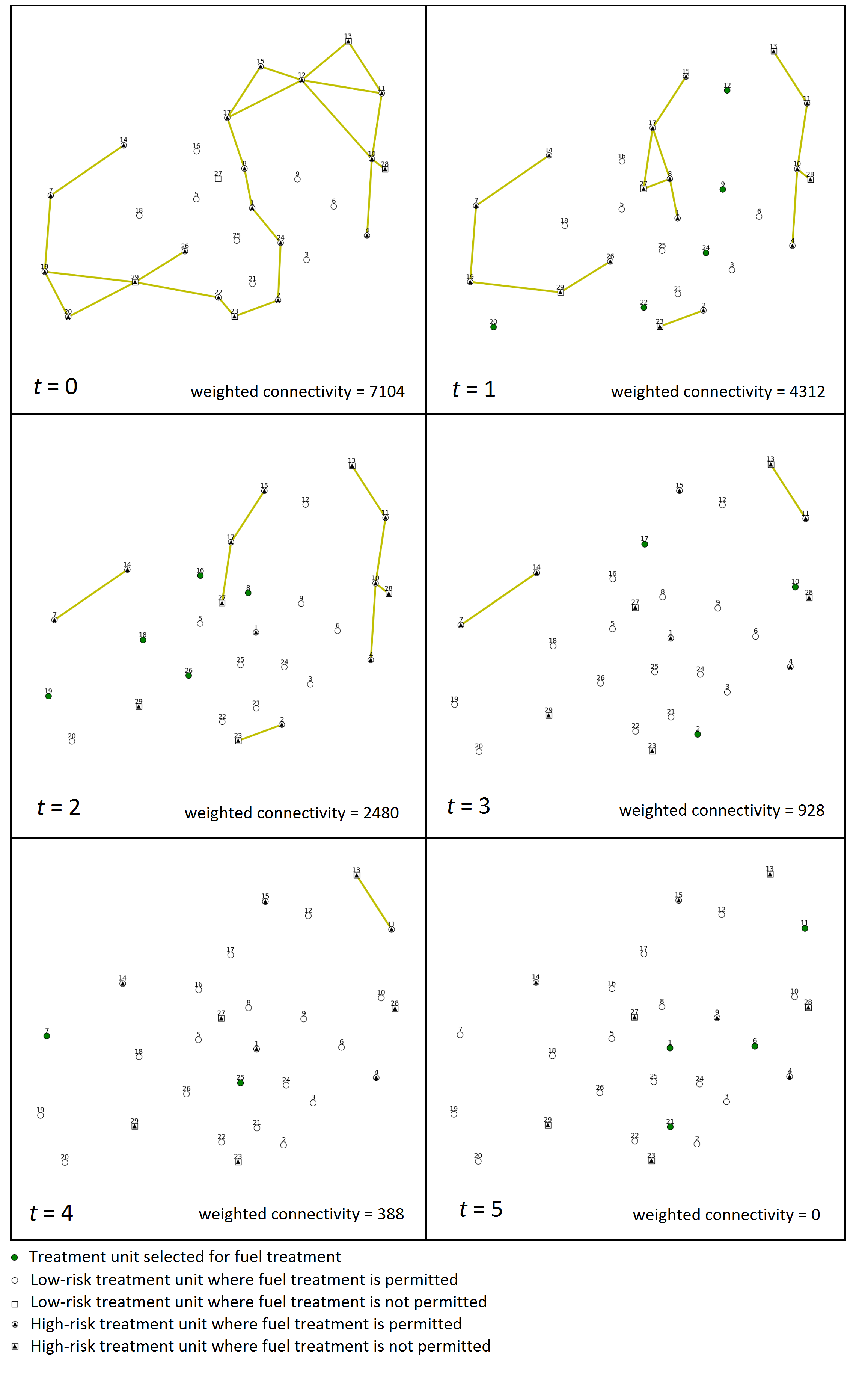}
\end{figure}

\begin{figure}
\protect\caption{The sequence of maps representing the fuel treatment schedule (in
years) for the test landscape\label{fig:A-map-represents-29-Fuel_scheduling}}

\includegraphics[scale=0.25]{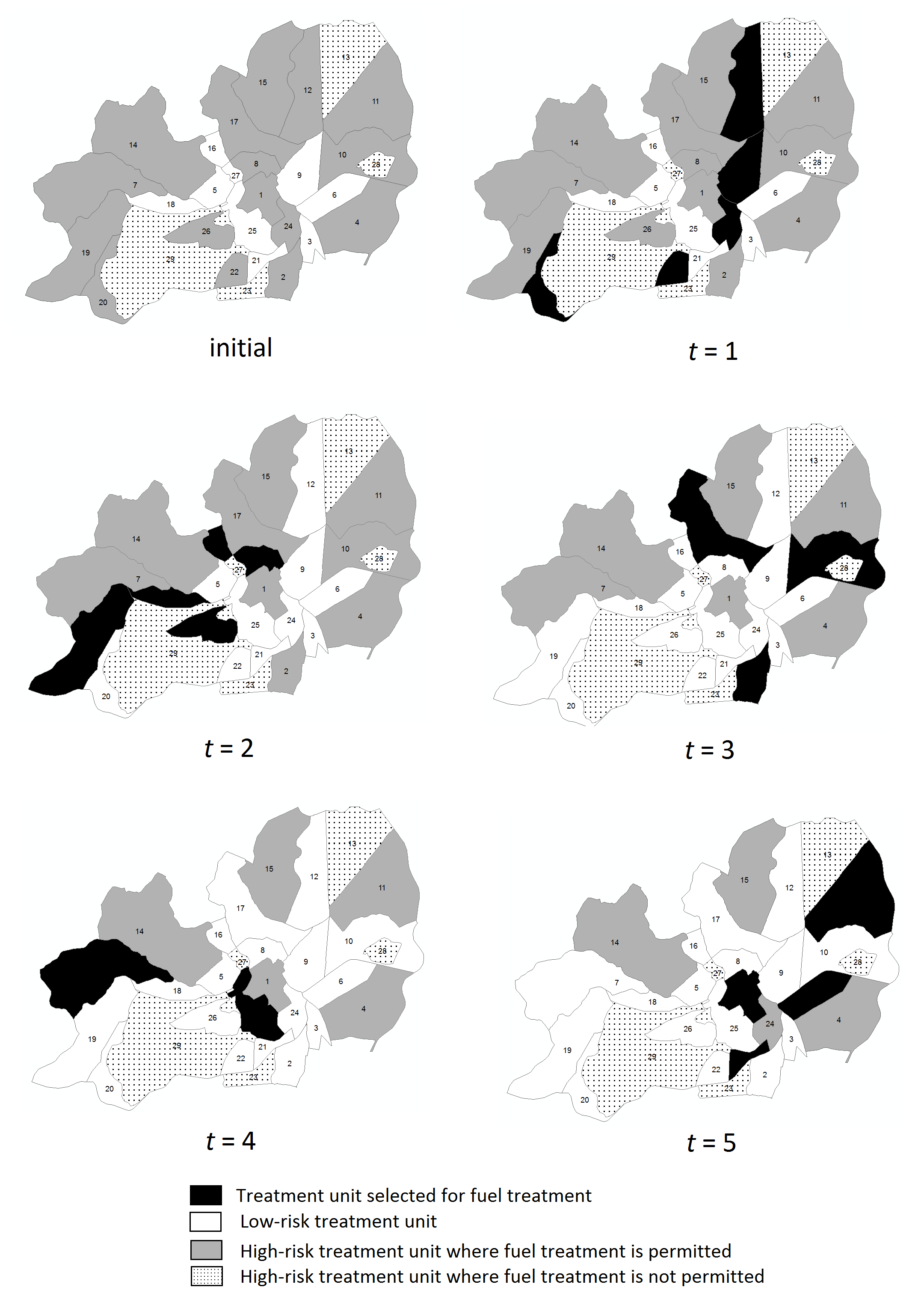}
\end{figure}

We evaluate the test landscape based on the data from Table \ref{tab:29-treatment-units}.
The rule is that if more than 50 percent of the treatment unit are
high-risk patches, then we consider it as a high-risk treatment unit.
Figure \eqref{fig:Fuel-treatment-schedule-29} and \eqref{fig:A-map-represents-29-Fuel_scheduling}
show the network and the related map representing the fuel treatment
schedule with 15 percent treatment level, starting from the \emph{t}
= 0 which represents the initial condition of the landscape. We can
treat the surrounding treatment units to break the connectivity of
high-risk units. When the patch within a treatment unit has reached
the maximum TFI, and no patch is below the minimum TFI, the treatment
units should be treated. This ecological requirement applies even
for the treatment units that do not contribute to the connectivity
of high-risk areas.

\section{An Australian case study }

\begin{figure}
\protect\caption{\label{fig:Location-of-the-case-study}}

\centering{}\subfloat[Location of the case study in the Barwon-Otway district of Victoria,
Australia\label{fig:Location-of-the-711map}]{

\protect\includegraphics[scale=0.2]{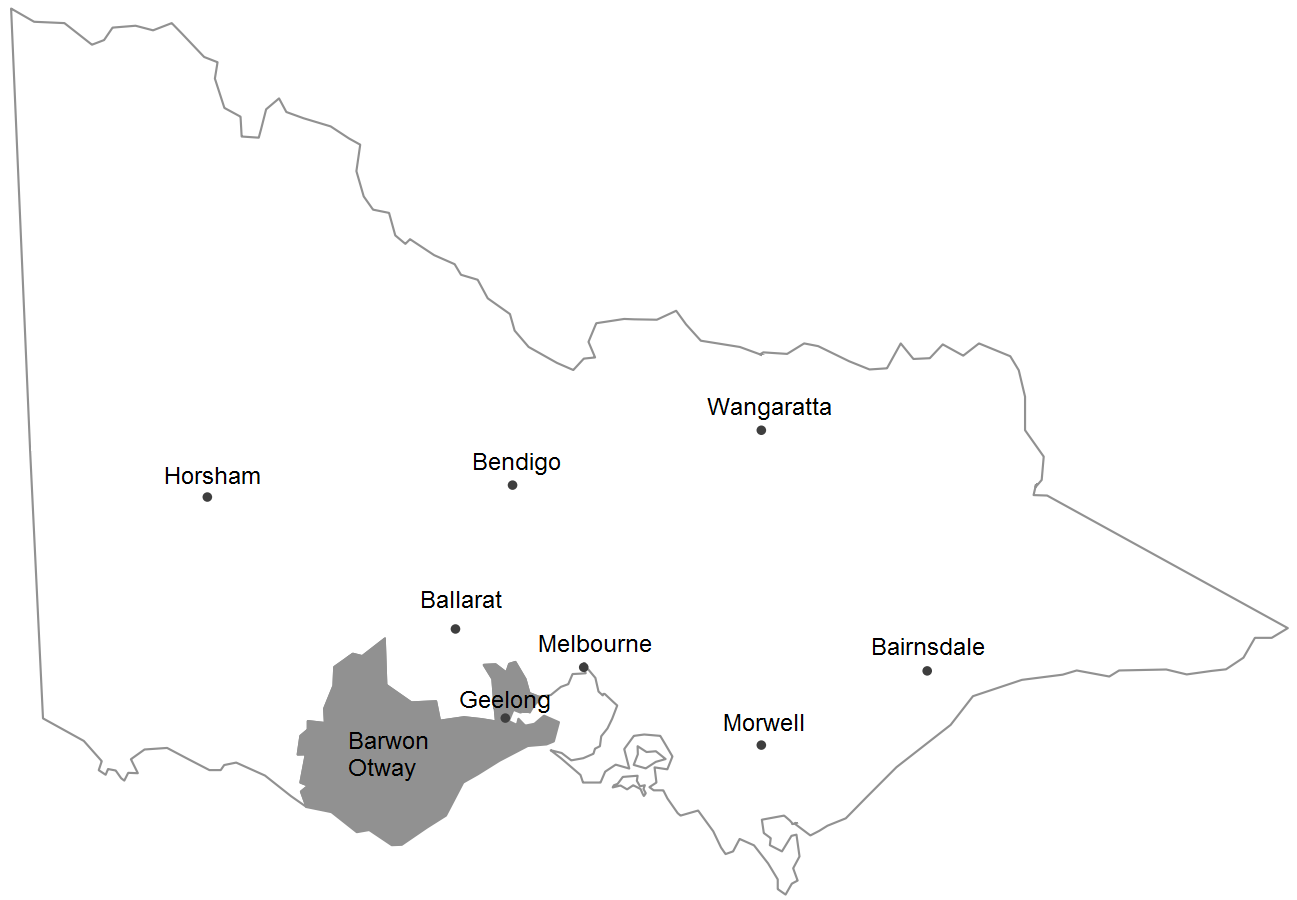}}\quad{}\subfloat[Map showing the distribution of the candidate treatment units within
case study area in the Barwon-Otway district of Victoria, Australia\label{fig:Map-showing-the-candidate-711}]{

\protect\includegraphics[scale=0.2]{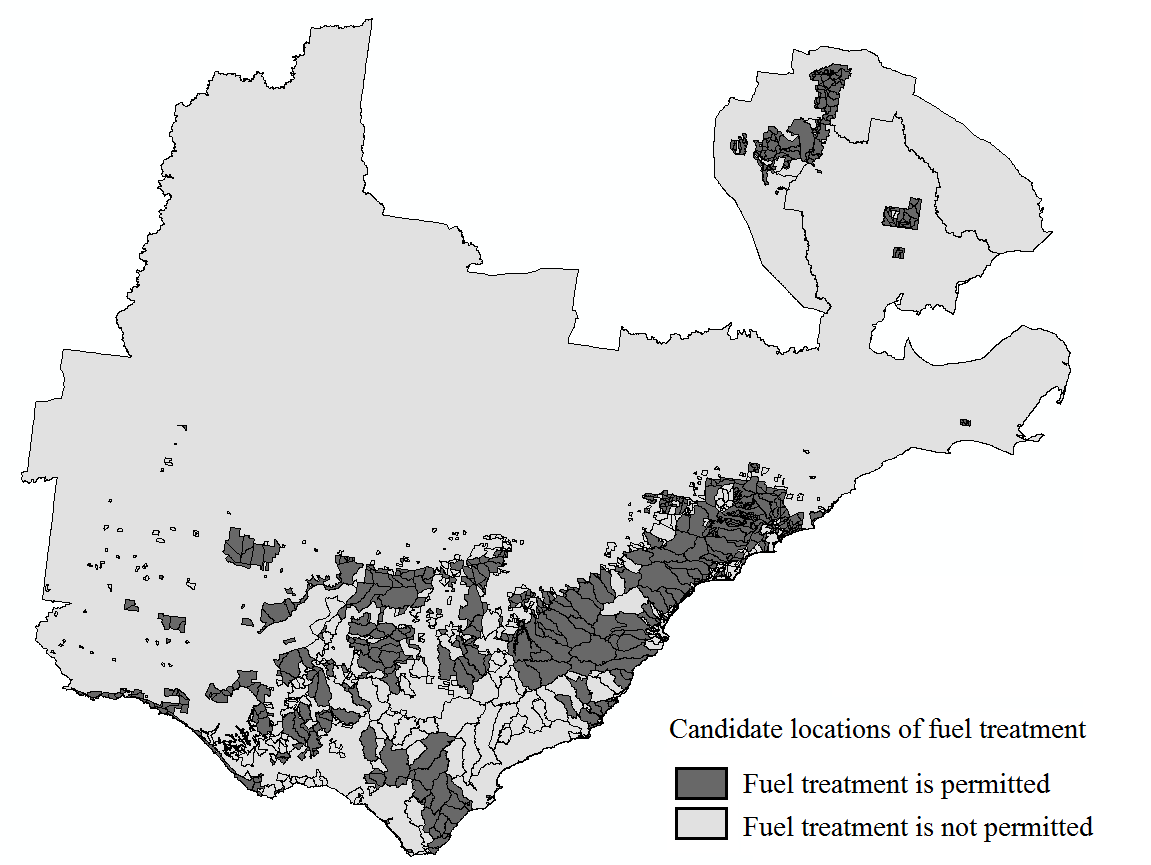}}
\end{figure}

\begin{figure}
\protect\caption{Solution of Phase 2: Maps showing the location of fuel treatment and
the spatial distribution of high-risk treatment units over time (in
years)\label{fig:Solution-of-Phase-711}}

\includegraphics[scale=0.22]{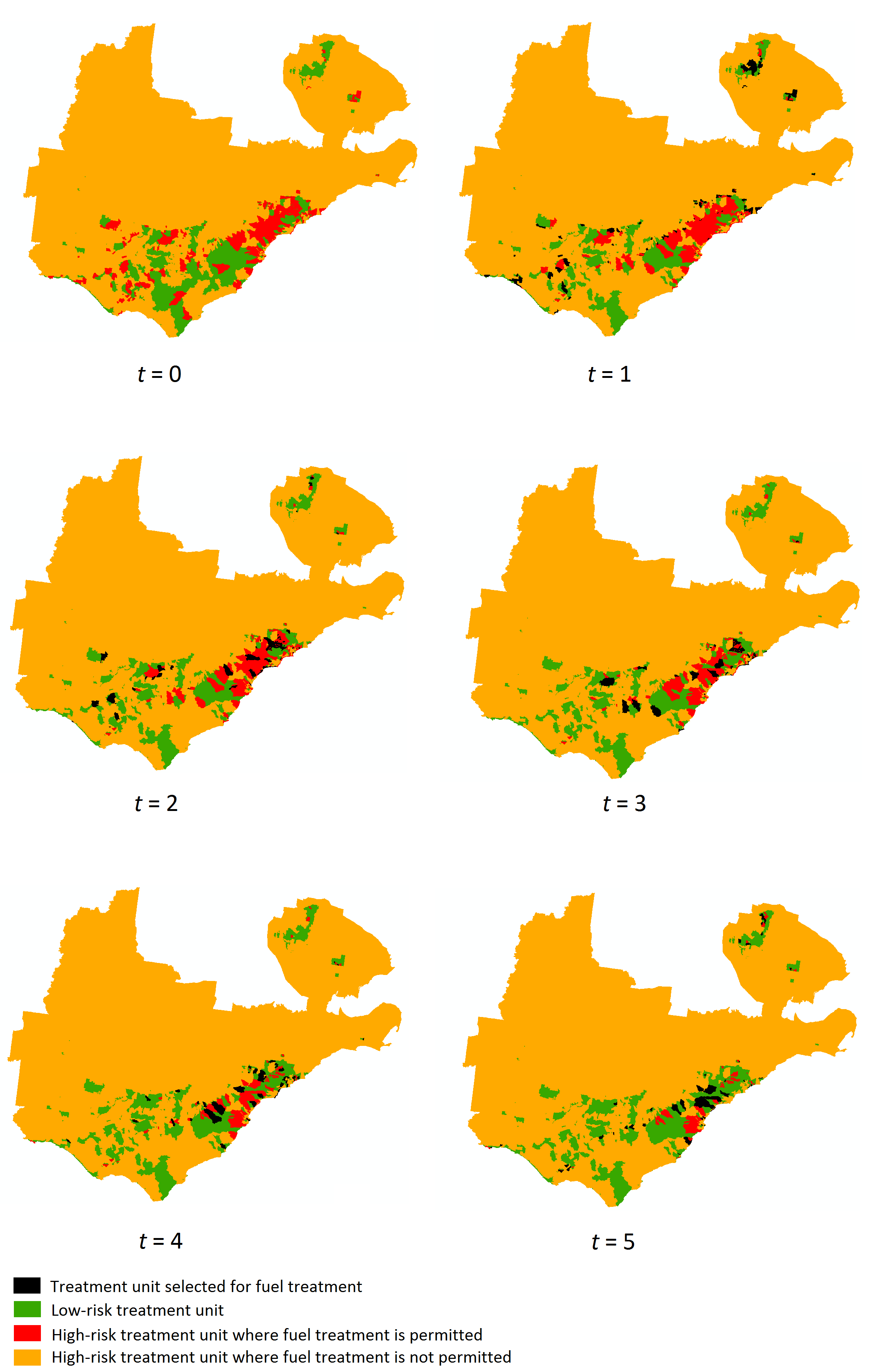}
\end{figure}

\begin{table}
\protect\caption{Ecological Vegetation Class (EVC) and the associated threshold age,
the minimum and the maximum TFI\label{tab:EVC-and-fuel-1}}

\bigskip{}

\centering{}%
\begin{tabular}{|l|>{\centering}m{0.7cm}|>{\centering}m{0.7cm}|>{\centering}m{1cm}|}
\hline 
{\footnotesize{}EVC name} & {\footnotesize{}min TFI}{\footnotesize \par}

{\footnotesize{}(year)} & {\footnotesize{}max TFI (year) } & {\footnotesize{}threshold (year)}\tabularnewline
\hline 
{\footnotesize{}Creekline Grassy Woodland} & {\footnotesize{}20 } & {\footnotesize{}150 } & {\footnotesize{}20}\tabularnewline
\hline 
{\footnotesize{}Hills Herb-rich Woodland} & {\footnotesize{}15 } & {\footnotesize{}150 } & {\footnotesize{}17}\tabularnewline
\hline 
{\footnotesize{}Creekline Herb-rich Woodland} & {\footnotesize{}15 } & {\footnotesize{}150 } & {\footnotesize{}17}\tabularnewline
\hline 
{\footnotesize{}Grassy Woodland} & {\footnotesize{}5} & {\footnotesize{}45} & {\footnotesize{}17}\tabularnewline
\hline 
{\footnotesize{}Valley Slopes Dry Forest} & {\footnotesize{}10} & {\footnotesize{}100} & {\footnotesize{}17}\tabularnewline
\hline 
{\footnotesize{}Sedgy Riparian Woodland} & {\footnotesize{}20} & {\footnotesize{}85} & {\footnotesize{}20}\tabularnewline
\hline 
{\footnotesize{}Scoria Cone Woodland} & {\footnotesize{}4} & {\footnotesize{}15} & {\footnotesize{}15}\tabularnewline
\hline 
{\footnotesize{}Wet Forest} & {\footnotesize{}45} & {\footnotesize{}300} & {\footnotesize{}45}\tabularnewline
\hline 
{\footnotesize{}Shrubby Wet Forest} & {\footnotesize{}25} & {\footnotesize{}150} & {\footnotesize{}25}\tabularnewline
\hline 
{\footnotesize{}Riparian Forest} & {\footnotesize{}10} & {\footnotesize{}80} & {\footnotesize{}22}\tabularnewline
\hline 
{\footnotesize{}Swampy Riparian Woodland} & {\footnotesize{}15} & {\footnotesize{}125} & {\footnotesize{}22}\tabularnewline
\hline 
{\footnotesize{}Riparian Scrub or Swampy Riparian Woodland Complex} & {\footnotesize{}10} & {\footnotesize{}80} & {\footnotesize{}16}\tabularnewline
\hline 
{\footnotesize{}Wet Sands Thicket} & {\footnotesize{}15} & {\footnotesize{}90} & {\footnotesize{}16}\tabularnewline
\hline 
{\footnotesize{}Stream Bank Shrubland} & {\footnotesize{}15} & {\footnotesize{}90} & {\footnotesize{}16}\tabularnewline
\hline 
{\footnotesize{}Cool Temperate Rainforest} & {\footnotesize{}45} & {\footnotesize{}999} & {\footnotesize{}45}\tabularnewline
\hline 
{\footnotesize{}Wet Heathland} & {\footnotesize{}12} & {\footnotesize{}45} & {\footnotesize{}12}\tabularnewline
\hline 
{\footnotesize{}Damp Heath Scrub} & {\footnotesize{}10} & {\footnotesize{}90} & {\footnotesize{}10}\tabularnewline
\hline 
{\footnotesize{}Damp Heath Scrub/Heathy Woodland Complex} & {\footnotesize{}10} & {\footnotesize{}90} & {\footnotesize{}10}\tabularnewline
\hline 
{\footnotesize{}Sand Heathland} & {\footnotesize{}8} & {\footnotesize{}45} & {\footnotesize{}8}\tabularnewline
\hline 
{\footnotesize{}Clay Heathland} & {\footnotesize{}10} & {\footnotesize{}45} & {\footnotesize{}10}\tabularnewline
\hline 
{\footnotesize{}Coastal Dune Scrub or Coastal Dune Grassland Mosaic} & {\footnotesize{}10} & {\footnotesize{}90} & {\footnotesize{}17}\tabularnewline
\hline 
{\footnotesize{}Coastal Headland Scrub} & {\footnotesize{}8} & {\footnotesize{}90} & {\footnotesize{}17}\tabularnewline
\hline 
{\footnotesize{}Coastal Headland Scrub/Coastal Tussock Grassland Mosaic} & {\footnotesize{}8} & {\footnotesize{}90} & {\footnotesize{}17}\tabularnewline
\hline 
{\footnotesize{}Coast Gully Thicket} & {\footnotesize{}10} & {\footnotesize{}90} & {\footnotesize{}17}\tabularnewline
\hline 
{\footnotesize{}Coastal Alkaline Scrub} & {\footnotesize{}10} & {\footnotesize{}70} & {\footnotesize{}17}\tabularnewline
\hline 
{\footnotesize{}Coastal Saltmarsh/Mangrove Shrubland Mosaic} & {\footnotesize{}8} & {\footnotesize{}90} & {\footnotesize{}14}\tabularnewline
\hline 
{\footnotesize{}Coastal Tussock Grassland} & {\footnotesize{}5} & {\footnotesize{}40} & {\footnotesize{}6}\tabularnewline
\hline 
{\footnotesize{}Heathy Woodland} & {\footnotesize{}5} & {\footnotesize{}45} & {\footnotesize{}35}\tabularnewline
\hline 
{\footnotesize{}Shrubby Woodland} & {\footnotesize{}10} & {\footnotesize{}45} & {\footnotesize{}35}\tabularnewline
\hline 
{\footnotesize{}Lowland Forest } & {\footnotesize{}8} & {\footnotesize{}80} & {\footnotesize{}20}\tabularnewline
\hline 
{\footnotesize{}Heathy Dry Forest } & {\footnotesize{}10} & {\footnotesize{}45} & {\footnotesize{}20}\tabularnewline
\hline 
{\footnotesize{}Shrubby Dry Forest } & {\footnotesize{}5} & {\footnotesize{}45} & {\footnotesize{}20}\tabularnewline
\hline 
{\footnotesize{}Grassy Dry Forest } & {\footnotesize{}5} & {\footnotesize{}45} & {\footnotesize{}15}\tabularnewline
\hline 
{\footnotesize{}Herb rich Foothill Forest } & {\footnotesize{}8} & {\footnotesize{}90} & {\footnotesize{}15}\tabularnewline
\hline 
{\footnotesize{}Shrubby Foothill Forest} & {\footnotesize{}8} & {\footnotesize{}90} & {\footnotesize{}15}\tabularnewline
\hline 
{\footnotesize{}Herb-rich Foothill Forest/Shrubby Foothill Forest
Complex} & {\footnotesize{}8} & {\footnotesize{}90} & {\footnotesize{}15}\tabularnewline
\hline 
{\footnotesize{}Damp Sands Herb Rich Woodland} & {\footnotesize{}10} & {\footnotesize{}90} & {\footnotesize{}17}\tabularnewline
\hline 
{\footnotesize{}Valley Grassy Forest} & {\footnotesize{}10} & {\footnotesize{}100} & {\footnotesize{}17}\tabularnewline
\hline 
{\footnotesize{}Plains Grassy Woodland} & {\footnotesize{}4} & {\footnotesize{}15} & {\footnotesize{}15}\tabularnewline
\hline 
{\footnotesize{}Alluvial Terraces Herb-Rich Woodland} & {\footnotesize{}4} & {\footnotesize{}15} & {\footnotesize{}15}\tabularnewline
\hline 
\end{tabular}
\end{table}

\begin{figure}
\protect\caption{The number of connected high-risk treatment units over time \label{fig:The-number-of-connection-1}}

\centering{}\includegraphics[scale=0.45]{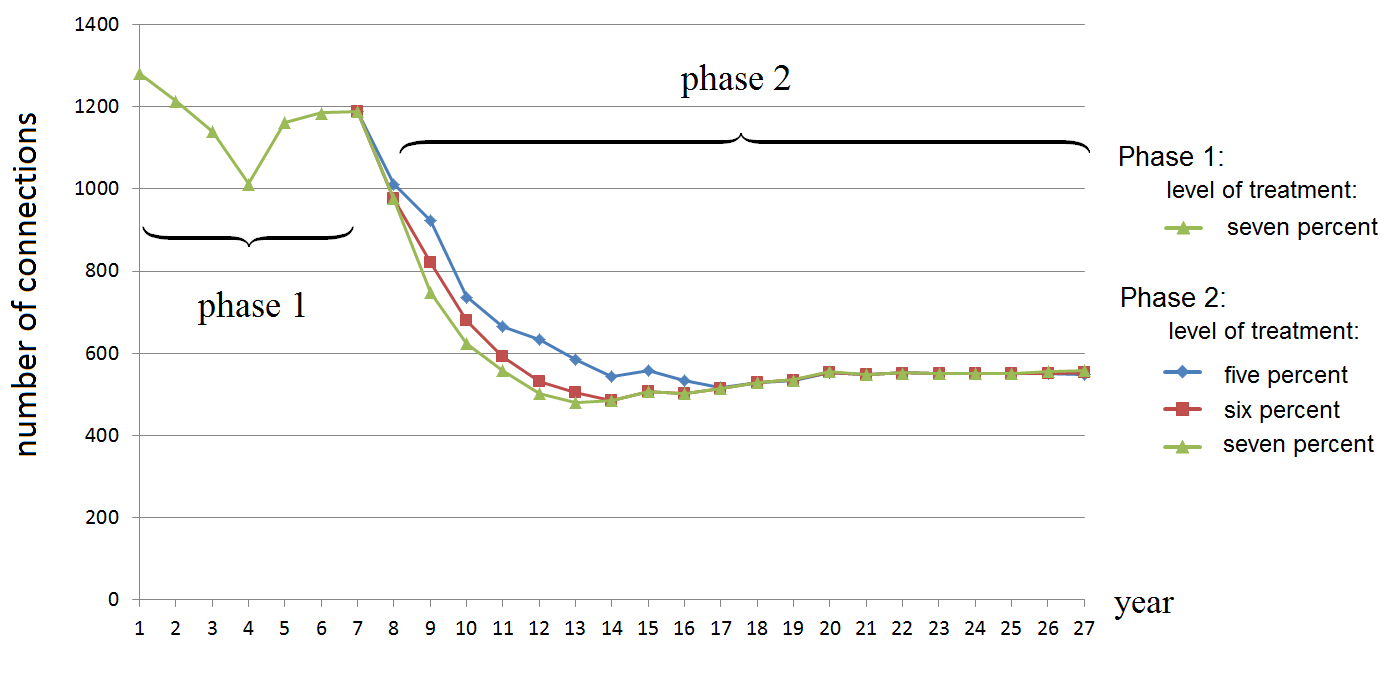}
\end{figure}

\begin{figure}
\protect\caption{The objective function values over time \label{fig:The-number-of-connection-1-1}}

\centering{}\includegraphics[scale=0.45]{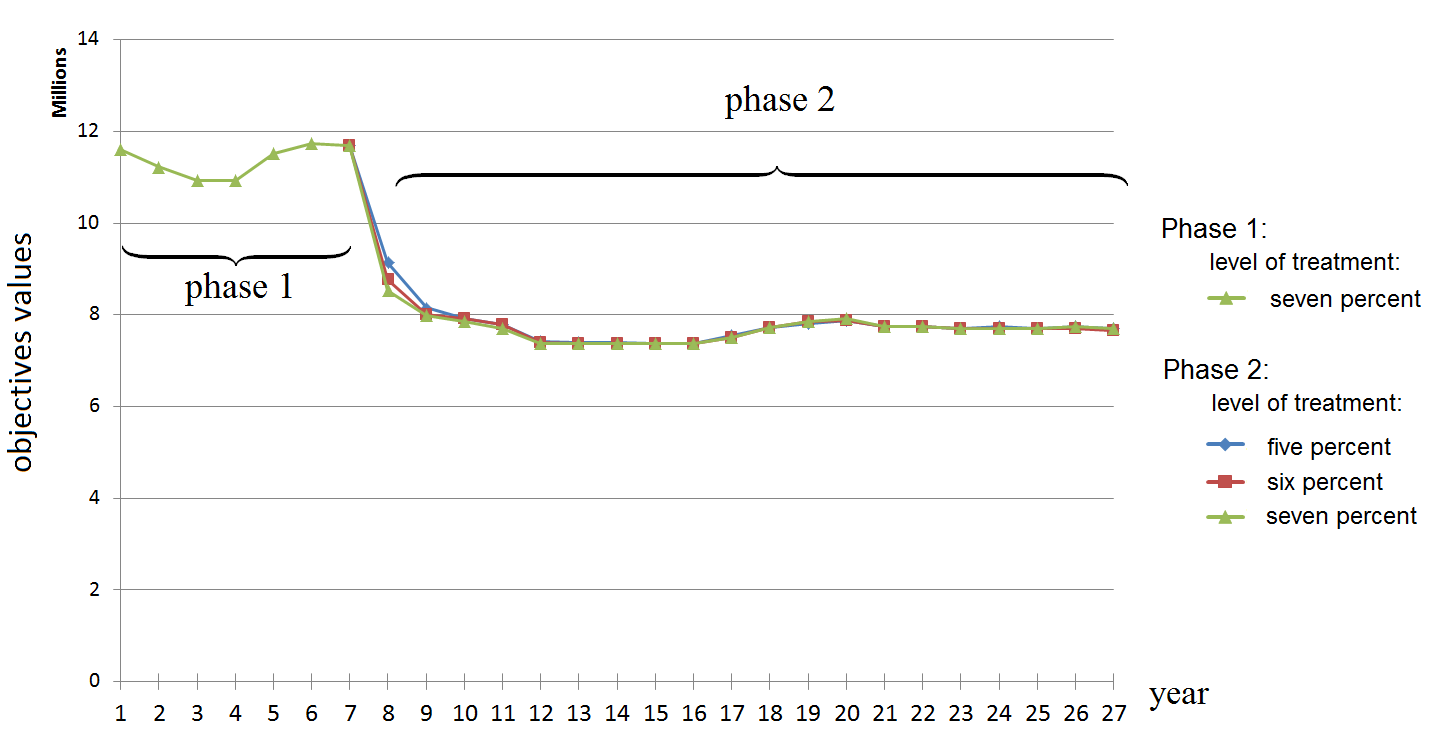}
\end{figure}

\begin{table}
\protect\caption{Computational comparison between the five, six and seven percent treatment
levels\label{tab:Computational-comparison-=002013}}

\bigskip{}

\centering{}%
\begin{tabular}{|>{\centering}p{3cm}|>{\centering}p{3cm}|>{\centering}p{3cm}|>{\centering}p{3cm}|}
\hline 
\multirow{2}{3cm}{%
Length of planning horizon%
} & \multicolumn{3}{c|}{%
Solution time (seconds) %
}\tabularnewline
\cline{2-4} 
 & %
five percent%
 & %
six percent%
 & %
seven percent%
\tabularnewline
\hline 
5 years%
 & %
22.32%
 & %
13.12%
 & %
11.72%
\tabularnewline
\hline 
10 years%
 & %
462.44%
 & %
38.29%
 & %
17.62%
\tabularnewline
\hline 
15 years%
 & %
4904.10%
 & %
752.11%
 & %
366.71%
\tabularnewline
\hline 
20 years%
 & %
26652.91%
 & %
9464.17%
 & %
2384.15%
\tabularnewline
\hline 
\end{tabular}
\end{table}

In this section, we apply the model discussed in Section \ref{sec:Model-formulation}
to an Australian case study. We use a real landscape with randomised
data containing treatable patches, grouped into 1197 treatment units.
Figure \ref{fig:Location-of-the-711map} illustrates the location
of the case study in the Barwon-Otway district of Victoria, Australia.
In this case study, we assume that we can only treat the public treatment
units. Figure \ref{fig:Map-showing-the-candidate-711} represents
the 711 candidate locations for fuel treatment. The data includes
area, vegetation type and age. The minimum TFI, maximum TFI and the
high-risk age threshold for each vegetation is summarised in Table
\ref{tab:EVC-and-fuel-1}. The vegetation types that do not pose any
threat such as aquatic vegetation types are excluded in this paper.
Threshold values are set to their assumed values to demonstrate our
approach rather than to provide an actual way of determining these
values. 

A set of connected treatment units is defined as a treatment unit
directly adjacent to another treatment unit, in other words, having
a shared boundary. It is acknowledged it is possible for treatment
units that are geographically separated to still be considered \textquoteleft connected\textquoteright{}
as a result of the spotting behaviour of particular bark fuel types
under given weather conditions. The provision of information regarding
bark fuel types and prevailing weather conditions for the case study
area would be a simple addition to model.

From the initial data, it was identified that 31 percent of the total
treatable area in the landscape is high-risk treatable treatment units
containing the patches that are over maximum TFI and no young patches.
Phase 1 is run for seven percent treatment level, and would need seven
years to achieve less than five percent high-risk treatment units
containing old patches in the landscape. In Phase 2, we run the model
presented in Section \ref{sec:Model-formulation} for five, six and
seven percent treatment levels. The solutions representing the high-risk
area over time and the location selected for fuel treatments each
year with seven percent treatment level can be seen in Figure \ref{fig:Solution-of-Phase-711}.
In this case study, we use the area of the two connected high-risk
treatment unit as a weight to determine the relative importance of
the connectivity. However, this weight can be determined in another
way, for example, by the proportion of the shared boundary between
two adjacent treatment units to the perimeter of the treatment units.
It can even be adjusted subjectively by the land manager if required.
Figure \ref{fig:The-number-of-connection-1} and \ref{fig:The-number-of-connection-1-1}
show that the number of connectivity of high-risk treatment units
in the landscape and the objective function values have the trends
to decrease over time. The average of the number of connections for
five, six and seven percent treatment levels are 608, 579 and 569,
respectively.

The model was solved using ILOG CPLEX 12.6 with the Python 2.7 programming
language using PuLP modeler. Computational experiments were performed
on Trifid, a V3 Alliance high-performance computer cluster. The computational
experiment used a single node with 16 cores of Intel Xeon E5-2670
64 GB of RAM. The comparison of computational time between the three
different treatment levels can be seen in Table \ref{tab:Computational-comparison-=002013}.
For the ten-year planning horizon, the computational time for the
three treatment levels is less than 15 minutes. For the longer planning
horizon, the computational time becomes longer. The optimal solution
can be obtained up to 20-year planning horizon.

\section{Conclusion}

In this paper, we have presented a mixed integer programming based
approach to schedule fuel treatments. The model determines when, and
where, to conduct the fuel treatment to reduce the fuel hazards in
the landscape whilst still meeting ecological requirements. The ecological
requirements considered in this paper are the minimum and maximum
Tolerable Fire Intervals (TFI) for the vegetation present. The model
includes multiple vegetation types and ages in the landscape and tracks
the age of vegetation in each treatment unit. To avoid deadlocks,
the rules that are applied in the model are either: the treatment
unit must be treated if there is an old patch in a treatment unit,
or the treatment unit cannot be treated if there is a young patch
in a treatment unit. In this study, spatial and temporal changes that
include multiple vegetation types in a more realistic polygon-based
network representation of the landscape are considered. These improve
upon previous work which was limited to a single vegetation type in
a regular grid, and create a more realistic approach to fuel treatment
planning for land managers. 

The model was illustrated in fuel treatment planning using real landscape
data from the Barwon-Otway district in south-west Victoria, Australia.
We ran the model for a 20-year planning horizon with five, six and
seven treatment levels. The total connectivity of high-risk regions
resulting from the three different treatment levels in the landscape
differs substantially for the first five years and differs slightly
after five years. Based on our experiments, using seven percent treatment
level, the high-risk regions in the landscape can be fragmented more
quickly than that of five and six percent, as expected. From the case
study, the solution of this complex multi-period model can be obtained
in a reasonable computational time (eight hours). Future work is planned
to extend this study by including habitat connectivity for fauna in
the landscape while still fragmenting high-risk areas.

\section*{\textcolor{black}{Acknowledgment}}

\textcolor{black}{The first author is supported by the Indonesian
Directorate General of Higher Education (1587/E4.4/K/2012). The second
author is supported by the Australian Research Council under the Discovery
Projects funding scheme (project DP140104246).}

\bibliographystyle{chicago}
\nocite{*}
\bibliography{RORH_arXiv}

\begin{thebibliography}{}

\bibitem[\protect\citeauthoryear{Ager, Vaillant, and Finney}{Ager
  et~al.}{2010}]{Ager2010_FEM}
Ager, A., N.~Vaillant, and M.~Finney (2010).
\newblock A comparison of landscape fuel treatment strategies to mitigate
  wildland fire risk in the urban interface and preserve old forest structure.
\newblock {\em Forest Ecology and Management\/}~{\em 259\/}(8), 1556--1570.

\bibitem[\protect\citeauthoryear{Boer, Sadler, Wittkuhn, McCaw, and
  Grierson}{Boer et~al.}{2009}]{Boer2009132}
Boer, M., R.~Sadler, R.~Wittkuhn, L.~McCaw, and P.~Grierson (2009).
\newblock Long-term impacts of prescribed burning on regional extent and
  incidence of wildfires-evidence from 50 years of active fire management in
  {SW} {A}ustralian forests.
\newblock {\em Forest Ecology and Management\/}~{\em 259\/}(1), 132--142.

\bibitem[\protect\citeauthoryear{Bradstock, Cary, Davies, Lindenmayer, Price,
  and Williams}{Bradstock et~al.}{2012}]{Bradstock2012Wildfires}
Bradstock, R., G.~Cary, I.~Davies, D.~Lindenmayer, O.~Price, and R.~Williams
  (2012).
\newblock Wildfires, fuel treatment and risk mitigation in australian eucalypt
  forests: Insights from landscape-scale simulation.
\newblock {\em Journal of Environmental Management\/}~{\em 105}, 66--75.

\bibitem[\protect\citeauthoryear{Burrows}{Burrows}{2008}]{Burrows20082394}
Burrows, N. (2008).
\newblock Linking fire ecology and fire management in {S}outh-{W}est
  {A}ustralian forest landscapes.
\newblock {\em Forest Ecology and Management\/}~{\em 255\/}(7), 2394 -- 2406.

\bibitem[\protect\citeauthoryear{Burrows and Wardell-Johnson}{Burrows and
  Wardell-Johnson}{2003}]{burrows2003fire}
Burrows, N. and G.~Wardell-Johnson (2003).
\newblock Fire and plant interactions in forested ecosystems of south-west
  {W}estern {A}ustralia.
\newblock {\em Fire in ecosystems of south-west {W}estern {A}ustralia: impacts
  and management\/}~{\em 2}, 225.

\bibitem[\protect\citeauthoryear{Cheal}{Cheal}{2010}]{Cheal2010}
Cheal, D. (2010).
\newblock Growth stages and tolerable fire intervals for {V}ictoria's native
  vegetation data sets.
\newblock In {\em Fire and adaptive management report no. 84}. Department of
  Sustainability and Environment, East Melbourne, Victoria, Australia.

\bibitem[\protect\citeauthoryear{Chung}{Chung}{2015}]{chung2015optimizing}
Chung, W. (2015).
\newblock Optimizing fuel treatments to reduce wildland fire risk.
\newblock {\em Current Forestry Reports\/}, 1--8.

\bibitem[\protect\citeauthoryear{Collins, Stephens, Moghaddas, and
  Battles}{Collins et~al.}{2010}]{Collins201024}
Collins, B., S.~Stephens, J.~Moghaddas, and J.~Battles (2010).
\newblock Challenges and approaches in planning fuel treatments across
  fire-excluded forested landscapes.
\newblock {\em Journal of Forestry\/}~{\em 108\/}(1), 24--31.

\bibitem[\protect\citeauthoryear{Ferreira, Constantino, and Borges}{Ferreira
  et~al.}{2011}]{ferreira2011sto}
Ferreira, L., M.~Constantino, and J.~Borges (2011).
\newblock A stochastic approach to optimize maritime pine (pinus pinaster ait.)
  stand management scheduling under fire risk. {A}n application in {P}ortugal.
\newblock {\em Annals of Operations Research\/}, 1--19.

\bibitem[\protect\citeauthoryear{Finney}{Finney}{2007}]{Finney2007702}
Finney, M. (2007).
\newblock A computational method for optimising fuel treatment locations.
\newblock {\em International Journal of Wildland Fire\/}~{\em 16\/}(6),
  702--711.

\bibitem[\protect\citeauthoryear{Garcia-Gonzalo, Pukkala, and
  Borges}{Garcia-Gonzalo et~al.}{2011}]{garcia2011}
Garcia-Gonzalo, J., T.~Pukkala, and J.~Borges (2011).
\newblock Integrating fire risk in stand management scheduling. an application
  to {M}aritime pine stands in {P}ortugal.
\newblock {\em Annals of Operations Research\/}, 1--17.

\bibitem[\protect\citeauthoryear{Kim, Bettinger, and Finney}{Kim
  et~al.}{2009}]{Kim2009253}
Kim, Y.-H., P.~Bettinger, and M.~Finney (2009).
\newblock Spatial optimization of the pattern of fuel management activities and
  subsequent effects on simulated wildfires.
\newblock {\em European Journal of Operational Research\/}~{\em 197\/}(1),
  253--265.

\bibitem[\protect\citeauthoryear{King, Bradstock, Cary, Chapman, and
  Marsden-Smedley}{King et~al.}{2008}]{King2008421}
King, K., R.~Bradstock, G.~Cary, J.~Chapman, and J.~Marsden-Smedley (2008).
\newblock The relative importance of fine-scale fuel mosaics on reducing fire
  risk in {S}outh-{W}est {T}asmania, {A}ustralia.
\newblock {\em International Journal of Wildland Fire\/}~{\em 17\/}(3),
  421--430.

\bibitem[\protect\citeauthoryear{Krivtsov, Vigy, Legg, Curt, Rigolot, Lecomte,
  Jappiot, Lampin-Maillet, Fernandes, and Pezzatti}{Krivtsov
  et~al.}{2009}]{Krivtsov20092915}
Krivtsov, V., O.~Vigy, C.~Legg, T.~Curt, E.~Rigolot, I.~Lecomte, M.~Jappiot,
  C.~Lampin-Maillet, P.~Fernandes, and G.~Pezzatti (2009).
\newblock Fuel modelling in terrestrial ecosystems: An overview in the context
  of the development of an object-orientated database for wild fire analysis.
\newblock {\em Ecological Modelling\/}~{\em 220\/}(21), 2915--2926.

\bibitem[\protect\citeauthoryear{Loehle}{Loehle}{2004}]{Loehle2004261}
Loehle, C. (2004).
\newblock Applying landscape principles to fire hazard reduction.
\newblock {\em Forest Ecology and Management\/}~{\em 198\/}(1-3), 261--267.

\bibitem[\protect\citeauthoryear{McCaw}{McCaw}{2013}]{McCaw2013217}
McCaw, L. (2013).
\newblock Managing forest fuels using prescribed fire - a perspective from
  southern {A}ustralia.
\newblock {\em Forest Ecology and Management\/}~{\em 294}, 217--224.

\bibitem[\protect\citeauthoryear{Minas, Hearne, and Martell}{Minas
  et~al.}{2014}]{Minas2014412}
Minas, J., J.~Hearne, and D.~Martell (2014).
\newblock A spatial optimisation model for multi-period landscape level fuel
  management to mitigate wildfire impacts.
\newblock {\em European Journal of Operational Research\/}~{\em 232\/}(2),
  412--422.

\bibitem[\protect\citeauthoryear{Penman, Christie, Andersen, Bradstock, Cary,
  Henderson, Price, Tran, Wardle, Williams, et~al.}{Penman
  et~al.}{2011}]{penman2011prescribed}
Penman, T., F.~Christie, A.~Andersen, R.~Bradstock, G.~Cary, M.~Henderson,
  O.~Price, C.~Tran, G.~Wardle, R.~Williams, et~al. (2011).
\newblock Prescribed burning: how can it work to conserve the things we value?
\newblock {\em International Journal of Wildland Fire\/}~{\em 20\/}(6),
  721--733.

\bibitem[\protect\citeauthoryear{Rachmawati, Ozlen, Reinke, and
  Hearne}{Rachmawati et~al.}{2015}]{Rachmawati2015}
Rachmawati, R., M.~Ozlen, K.~Reinke, and J.~Hearne (2015).
\newblock A model for solving the prescribed burn planning problem.
\newblock {\em SpringerPlus\/}~{\em 4\/}(1).

\bibitem[\protect\citeauthoryear{Rytwinski and Crowe}{Rytwinski and
  Crowe}{2010}]{Rytwinski2010}
Rytwinski, A. and K.~A. Crowe (2010).
\newblock A simulation-optimization model for selecting the location of
  fuel-breaks to minimize expected losses from forest fires.
\newblock {\em Forest Ecology and Management\/}~{\em 260\/}(1), 1 -- 11.

\bibitem[\protect\citeauthoryear{Schmidt, Taylor, and Skinner}{Schmidt
  et~al.}{2008}]{Schmidt20083170}
Schmidt, D., A.~Taylor, and C.~Skinner (2008).
\newblock The influence of fuels treatment and landscape arrangement on
  simulated fire behavior, {S}outhern {C}ascade range, {C}alifornia.
\newblock {\em Forest Ecology and Management\/}~{\em 255\/}(8-9), 3170--3184.

\bibitem[\protect\citeauthoryear{Van~Wagtendonk}{Van~Wagtendonk}{1995}]{vanwagtendonk1995large}
Van~Wagtendonk, J.~W. (1995).
\newblock Large fires in wilderness areas.
\newblock {\em United States Department of Agriculture Forest Service General
  Technical Report Int\/}, 113--116.

\bibitem[\protect\citeauthoryear{Wei and Long}{Wei and
  Long}{2014}]{wei2014schedule}
Wei, Y. and Y.~Long (2014).
\newblock Schedule fuel treatments to fragment high fire hazard fuel patches.
\newblock {\em Mathematical and Computational Forestry \& Natural-Resource
  Sciences (MCFNS)\/}~{\em 6\/}(1), 1--10.

\end{thebibliography}
%RORH_arXiv

\vspace{3cm}

\noindent Ramya Rachmawati

\noindent Mathematics Department Faculty of Mathematics and Natural
Sciences, University of Bengkulu, Bengkulu, Indonesia

\noindent School of Mathematical and Geospatial Sciences RMIT University,
Melbourne, Australia

\noindent (ramya.rachmawati@rmit.edu.au)

\noindent \medskip{}

\noindent Melih Ozlen

\noindent School of Mathematical and Geospatial Sciences RMIT University,
Melbourne, Australia

\noindent (melih.ozlen@rmit.edu.au)

\noindent \medskip{}

\noindent Karin J. Reinke

\noindent School of Mathematical and Geospatial Sciences RMIT University,
Melbourne, Australia

\noindent (karin.reinke@rmit.edu.au)

\noindent \medskip{}

\noindent John W. Hearne

\noindent School of Mathematical and Geospatial Sciences RMIT University,
Melbourne, Australia

\noindent (john.hearne@rmit.edu.au)
\end{document}